\documentclass[10pt,twoside]{amsart}

\title[The number of triangulations of a planar point set]
            {A better upper bound on the \\
         number of triangulations of a planar point set}
\author{Francisco Santos}
\author{               Raimund Seidel}

\address{Francisco Santos,
    Dept.~de Matem\'aticas,
    Facultad de Ciencias,
    Universidad de Cantabria,
    E-39005 Santander, Spain.}
\email{santos@matesco.unican.es.}
\thanks{The first author was partially supported by grant
    BFM2001--1153 of the Spanish  Direcci\'on General de Ense\~nanza Superior
    e Investigaci\'on Cient\'{\i}fica.}

\address{Raimund Seidel,
    FR 6.2 Informatik,
    Universit\"at des Saarlandes,
    Postfach 151150,
    D-66041 Saarbr\"ucken, Germany.} 
\email{rseidel@cs.uni-sb.de.}

\date{February 2002}

\font\bbb=msbm10 scaled 1200

\usepackage{latexsym}
\usepackage{epsf}

\newtheorem{theorem}{Theorem}

\newtheorem{remark}[theorem]{Remark}
\newtheorem{lemma}[theorem]{Lemma}

\newcommand{\bigO}{\operatorname{O}}
\newcommand{\reals}{\hbox{\bbb R}}
\newcommand{\nats}{\hbox{\bbb N}}
\newcommand{\conv}{\operatorname{conv}}

\begin{document}
\maketitle

\begin{abstract}

We show that a point set of cardinality $n$ in the plane cannot be the
vertex set of more than $59^{n} \bigO(n^{-6})$ 
straight-edge triangulations of its convex
hull. This improves the previous upper bound of $276.75^{n
+\bigO(\log(n))}$.
\end{abstract}

\section*{Introduction}

A triangulation of a finite point set 
$A\subset\reals^2$ in the Euclidean plane
is a geometric simplicial complex covering
$\conv(A)$ whose vertex set is precisely $A$. Equivalently, it is a maximal
non-crossing straight-edge graph with vertex set $A$.
In this paper we prove that a point set of cardinality $n$ 
cannot have more than $59^n\bigO(n^{-6})$
triangulations. 

An upper bound of type $2^{\bigO(n)}$ for this number
is a consequence of the general results of \cite{ACNS}. 
Upper bounds of
$173\,000^n$, $7\,187.52^n$ and $276.75^{n+ \bigO(\log(n))}$ 
have been given,
respectively, in \cite{Sm}, \cite{Se} and \cite{DS}.
The precise statement of our new upper bound is:

\begin{theorem}
\label{thm:main}
The number of triangulations of a planar point set is bounded above by
$$
\frac{59^v\cdot 7^b}{{v+b+6 \choose 6}},
$$
where $v$ and $b$
denote the numbers of \emph{interior} and \emph{boundary} points,
meaning by this points lying in the interior
and the boundary of the convex hull, respectively. 
\end{theorem}

In Section \ref{sec:main} we prove this result and in Section \ref{sec:review}
we briefly review what is known about the maximum and minimum number of
triangulations of point sets of fixed cardinality. 
In particular, we mention
that every point set in general position has at least
$\Omega(2.012^{n})$, as proved in \cite{AHN}.
As a reference, compare these upper and lower bounds to the number of
triangulations of $n$ points in convex position, which is 
the Catalan number $C_{n-2}=\frac{1}{n-1}{2n-4 \choose n-2}
= \Theta(4^n n^{-\frac{3}{2}})$.

\section{Proof of the upper bound}
\label{sec:main}
We will assume that our point set 
is in general position, i.e., that no three of
its points are collinear. This is no loss of generality because if $A$
is not in general position and we perturb it to a point set $A'$ in
general position without making boundary points go to the interior,
then every triangulation of $A$ is a triangulation of $A'$ as well. In
particular, the maximum number of triangulations of planar point sets
of given cardinality is achieved in general position. The same 
may not be true in higher dimension.

For the proof of Theorem~\ref{thm:main} we will need the fact that
the total number of vertices in a triangulation is bounded by a linear
combination of the number of low-degree vertices:

\begin{lemma}
\label{lemma:euler}
Let $T$ be a triangulation in the plane. For each integer $i\ge 3$, let
$v_i$ denote the number of interior vertices of degree $i$ in $T$. For
each integer $j\ge 2$, let $b_j$ denote the number of boundary
vertices of degree $j$ in $T$. Then
\[
\sum (6-i) v_i + \sum (4-j) b_j = 6\,, 
\]
and therefore
$$
|A| + 6 \leq 4 v_3 + 3 v_4 + 2 v_5 + v_6 + 3 b_2 + 2 b_3 + b_4\,.
$$
\end{lemma}

\begin{proof}
  Let $v$ and $b$ be the numbers of interior and boundary vertices in $T$,
  respectively. Let $e$ and $t$ be the numbers of edges and triangles in $T$.
  Counting the edges of $T$ according to their incidences to triangles shows
  that $3t = 2e - b$.  Euler's formula says that $t-e+b+v=1$. These two
  equations give:
$$
6v + 4b = 6 + 2e.
$$
On the other hand, counting the edges of
$T$ according to their incidences to vertices shows that
$$
2e = \sum i v_i + \sum j b_j. 
$$
Substituting this into the previous equality and noting that $v=\sum_{i\geq
  3} v_i$ and $b=\sum_{j\geq 2}b_j$ gives the first claimed equation.  Adding
$|A|=\sum_i v_i + \sum_j b_j$ on both sides of this equality and dropping the
negative summands on the left hand side yields the claimed inequality.

\end{proof}

Let $T$ be a triangulation of $A$ and $p\in A$ be one of its
points. We say that a triangulation $T'$ of $A\setminus\{p\}$ is
obtained by \emph{deleting $p$ from $T$} if all the edges of $T$ not
incident to $p$ appear in $T'$. In the same situation we also say
that $T$ is obtained by \emph{inserting $p$ into $T'$}. Observe that
neither the deletion nor the insertion of a point into a triangulation
is a unique process:
 more than one triangulation of $A'$ can result
from the deletion of $p$ from $T$ and more than one triangulation
of $A'$ can result from the insertion of $p$ into $T'$.
However, these numbers can be bounded in terms of Catalan numbers.
\begin{lemma}
\label{lemma:delete}
Let $T$ be a triangulation of $A$ and let $p$ be a vertex in $A$ with
degree $i$.  The number of triangulations of $A'=A\setminus \{p\}$ 
that can be obtained from $T$ by deleting $p$ is at least 1 and at most
$C_{i-2}=\frac{1}{i-1}{2i-4 \choose i-2}$.
\end{lemma}
\begin{proof}
It suffices to
note that the number in question is the number of ways in which the
area formed by the intersection of $\conv(A')$ and the triangles in $T$
incident to $p$ can be triangulated.
\end{proof}

\begin{lemma}
\label{lemma:insert}
Let $T'$ be a triangulation of $A\setminus \{p\}$. For each $i\in
\nats$, let $h_i$ be the number of triangulations of $A$ in which $p$
has degree $i$ and which can be obtained by inserting $p$ into $T'$.
\begin{itemize}
\item If $p$ is an interior point of $A$, then
$h_i \leq C_{i-1}-C_{i-2}=\frac{3}{2i-3}{2i-3 \choose i-3}$. In particular
$h_3 = 1$, $h_4 \le 3$,
$h_5 \le 9$ and $h_6 \le 28$.
\item If $p$ is a boundary point of $A$, then
$h_i \leq C_{i-2}=\frac{1}{i-1}{2i-4 \choose i-2}$. In particular
$h_2 \le 1$, $h_3 \le 1$ and
$h_4 \le  2$.
\end{itemize}

\end{lemma}

\begin{proof}
Let us first assume that $p$ is interior. Let $\Delta$ be the triangle of
$T'$ that contains $p$. After inserting 
$p$ with degree $i$ in $T'$, the union of the triangles incident to $p$ is
a starshaped polygon $Q$, with $p$ in its kernel, with no other
point of $A$ in its interior, and obtained as the union of $i-2$ triangles
from $T'$. Conversely, any polygon with those properties
provides a way of inserting $p$ with degree $i$.

Any triangle $t$ of $T'$ in such a polygon will be
\emph{visible} from $p$, meaning that $\conv(t\cup\{p\})$
contains no vertex of $T'$ in its interior.  Let $G'$ be the dual graph
of $T'$, whose nodes are the triangles in $T'$ with two of them adjacent
if the triangles share a common edge.  Let $W$ be the subgraph in $G'$
induced by the triangles visible from $p$.  This subgraph cannot contain
a cycle (by non-degeneracy $p$ is not collinear with any two vertices
of $T'$) and thus it is a forest.  Let $V$ be the tree in $W$ that contains
$\Delta$. The number $h_i$ coincides with the number of subtrees of $V$ with
$i-2$ nodes that include $\Delta$.

We can view $V$ as a planted
tree with root $\Delta$, which has degree at most $3$, and all of whose
subtrees are binary trees.  The number of $(i-2)$-node subtrees of $V$
that contain the root is upperbounded by the number of such subtrees
of the infinte tree $Z$ whose root has exactly three children each of which
is root of an infinte binary tree. In turn, this number equals the number of
binary trees with $(i-1)$ nodes whose right spine is not empty,
because $Z$ can be also described as the tree obtained 
contracting the first edge of the right spine in the infinite binary tree.
The number of those trees is clearly $C_{i-1}-C_{i-2}$.

Suppose now that $p$ is a boundary point. Let $e_1,\ldots,e_k$ be
the edges for which $p$ is beyond, in order along the boundary of $T'$, and
let $t_i=\conv(e_i\cup\{p\})$.  Enlarge the triangulation $T'$ by the
traingles $t_1,\ldots,t_k$ to a triangulation $T''$ and proceed
as in the first case but with $t_1$ playing the role of $\Delta$.
The desired upper bound then turns out to be given by the
number of $(i-1)$-node binary trees with exactly $k$ nodes on the
right spine.
\end{proof}

\begin{proof}[Proof of Theorem \ref{thm:main}:]
Let $N(v,b)$ denote, for every pair of integers
$v\ge 0$ and $b\ge 3$, the maximum number of triangulations among all
point sets with $v+b$ points and with \emph{at most} $v$ of them
interior.  We will prove by induction on $v+b$ that 
$N(v,b)\le {59^{v}\cdot
7^{b}}/{{v+b+6\choose 6}}$. Induction
starts with $b=3$ and $v=0$, which gives 
${59^{0}\cdot 7^{3}}/{{0+3+6\choose 6}} = 49/12\ge 1 = N(3,0)$.

Let $A$ be a point set with $b$ boundary points and $v$ interior points.
For each $i\ge 3$, let $V_i$ denote the sum over all triangulations of
$A$ of the numbers of interior vertices of degree $i$.  For each $j\ge
2$, let $B_j$ denote the sum over all triangulations of $A$ of the
numbers of boundary vertices of degree $j$. Let $N$ be the number of
triangulations of $A$.  Observe that deleting an interior point from
$A$ gives a point configuration with $b$ boundary points and $v-1$
interior points, while deleting a boundary point gives a point
configuration with $v+b-1$ points, \emph{at most $v$ of which} are
interior. 

The number of triangulations of $T$ in which a certain vertex
$p$ has degree $i$ is at most equal to the number of ways of inserting
$p$ with degree $i$ in triangulations of $A\setminus\{p\}$. The
inequality may by strict since insertions from different
triangulations of $A\setminus\{p\}$ can lead to the same triangulation
of $A$. Then, Lemma \ref{lemma:insert} implies that
$$
V_3 \le v N(v-1,b),
\qquad\qquad 
V_4 \le 3 v N(v-1,b),
$$
$$
V_5 \le 9 v N(v-1,b),
\qquad\qquad 
V_6 \le 28 v N(v-1,b),
$$
$$
B_2 \le b N(v,b-1), \qquad 
B_3 \le b N(v,b-1)
\qquad \hbox{and}\qquad 
B_4 \le 2 b N(v,b-1).
$$
On the other hand, from Lemma
\ref{lemma:euler} we have that 
$$
(6 + v + b) N 
  \le 4 V_3 + 3 V_4 + 2 V_5 + V_6 + 3 B_2 + 2 B_3 + B_4
  \le 59 v N(v-1,b) + 7 b N(v,b-1).
$$

By inductive hypothesis,
$N(v-1,b)\le \frac{59^{v-1}\cdot 7^{b}}{{v+b+5\choose 6}}$ and
$N(v,b-1)\le \frac{59^{v}\cdot 7^{b-1}}{{v+b+5\choose 6}}$. Hence,
$$
(6 + v + b) N \le (v+b) \frac{59^{v}\cdot 7^{b}}{{v+b+5\choose 6}}
=(6 + v + b) \frac{59^{v}\cdot 7^{b}}{{v+b+6\choose 6}}.
$$
\end{proof}

\begin{remark}\rm
\label{rem:not-all-points}
In the theory of secondary polytopes and the so-called Baues problem
(see, for example, \cite{BFS} and \cite{Re}) 
it is natural to consider as triangulations of $A$
those covering $\conv(A)$ and 
with vertex set \emph{contained} in $A$, allowing not to use all of
the interior
points as vertices. We can bound the number of triangulations
in this setting by adding the bounds of Theorem \ref{thm:main}
for the different subsets of interior vertices. This gives the
following upper bound, where $N(v,b)$ is as in the proof of Theorem
\ref{thm:main}:
$$
\sum_{i=0}^{v}{v \choose i} N(i,b) \le
\sum_{i=0}^{v}{v \choose i}\frac{59^{i}\cdot 7^{b}}{{i+b+6\choose 6}} \le
\sum_{i=0}^{v}{v \choose i}\frac{59^{i}\cdot 7^{b}}{{b+6\choose 6}} =
\frac{ 60^{v} \cdot 7^{b}}{{b+6 \choose 6}}.
$$
\end{remark}

\section{The number of triangulations in the plane. State of the art}
\label{sec:review}

Let $T(n)$ and $t(n)$ denote the maximum and minimum number of triangulations
among all point sets in the plane 
in general position and of cardinality $n$. 
For $t(n)$ it is necessary to assume general position
since $n-1$ points on a line produce a point set with
only 1 triangulation.
For $T(n)$,
general position is no loss of generality. 
The following table, taken from \cite{Aich}, 
gives $T(n)$ and $t(n)$ for $n=3,\dots,10$, compared to the number
of triangulations of the convex $n$-gon:

\begin{table}[h]
\begin{tabular}{|c|c|c|c|c|c|c|c|c|}
\hline
   $n$    & 3 & 4 & 5 &  6 &  7 &   8 &   9 &   10 \\ 
\hline
$t(n)$    & 1 & 1 & 2 &  4 & 11 &  30 &  89 &  250 \\ 
$C_{n-2}$ & 1 & 2 & 5 & 14 & 42 & 132 & 429 & 1430 \\ 
$T(n)$    & 1 & 2 & 5 & 14 & 42 & 150 & 780 & 4550 \\ 
\hline
\end{tabular}
\end{table}

Concerning the asymptotic behaviour of $t(n)$ and $T(n)$ we know that:
\begin{eqnarray}
\Omega(2.0129^n) &\le\quad t(n) \quad\le& \bigO(12^{n/2})=\bigO(3.46410^n),\\ 
\Omega(8^n n^{-7/2})&\le\quad T(n) \quad\le& \bigO(59^n n^{-6})
\end{eqnarray}
Compare this with $C_{n-2}= \Theta(4^n n^{-\frac{3}{2}})$ for the convex
$n$-gon.
The lower bound for $t(n)$ comes from \cite{AHN}. The upper bound for $T(n)$
is our Theorem \ref{thm:main}. The other two bounds come from the
computation of the number of triangulations of the following point
sets:
\begin{itemize}
\item {\bf A double chain:} Let $A$ consist of two convex chains 
of $k=n/2$ points each, facing one another and so that every pair of
segments in different chains are visible from one another.
See the center picture in Figure \ref{fig:extremal}, for the case $k=9$. The
edges drawn
in the figure are ``unavoidable'', i.e., present in every triangulation.
They divide $A$ into two convex $k$-gons, with
$C_{k-2}$ triangulations each, and a non-convex $2k$-gon which is easily seen
to have 
${2k-2 \choose k-1}$ triangulations (see \cite{GNT}).
Hence, the number of triangulations of $A$ is:
\[
{2k-2 \choose k-1} {C_{k-2}}^2 =\Theta(64^{k}k^{-\frac{7}{2}})=
\Theta(8^{n}n^{-\frac{7}{2}}).  
\]
\begin{figure}
\begin{center}
\epsfxsize=10 cm
\epsfbox{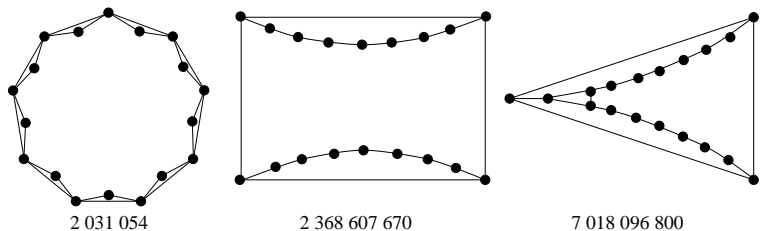}
\caption{The point sets with (asymptotically) the greatest and lowest number
  of triangulations known.
 The number of triangulations of each, for $n=18$, is shown.}
\label{fig:extremal}
\end{center}
\end{figure}

\item {\bf A double circle:} Let $A$ be a 
 convex $k$-gon ($k=n/2$)  together with 
an interior point sufficiently close to each boundary edge.
See the left part of Figure \ref{fig:extremal}, for the case $k=9$. 
Again, the edges drawn are unavoidable, and triangulating $A$ is the same as
triangulating the central non-convex $n$-gon. An inclusion-exclusion
argument (see Proposition 1 in \cite{HN97}) 
gives the exact number of
triangulations of this polygon, which is
\[
\sum_{i=0}^{k} (-1)^{i} {k \choose i} C_{2k-i-2}\le 12^k =12^\frac{n}{2}.
\]
%
\end{itemize}

It is interesting to observe that the double circle actually gives the minimum
possible number of triangulations for $n\le 10$. 
See \cite{Aich}, where this is conjectured to be true for all $n$.
(If $n$ is odd, the double circle has to be modified with
an extra interior point.) 

On the other hand, the double chain has only $6$, $80$ and $1750$ 
triangulations for $n=6$, $8$ and $10$ respectively,
which is less than $T(n)$. There is a simple way to modify it and
get more triangulations, as shown in the right picture of Figure
\ref{fig:extremal}. The big non-convex polygon is a
double chain with $n-2$ vertices. This modified double chain has exactly
\[
{2k-4 \choose k-2} {C_{k-1}}^2 =\frac{(2k-3)(2k-2)}{k^2}
{2k-2 \choose k-1} {C_{k-2}}^2 
\]
triangulations, where $n=2k$ as in the previous examples.
That number is (asymptotically) four times the number of
triangulations of the double chain. Still, the modified double chain has $8$,
$150$ and $3920$ triangulations for $n=6,8,10$, which is not (always) the
maximum. 
The numbers of triangulations of these configurations for $n=18$ appear in
Figure \ref{fig:extremal}. The greatest number of triangulations for $n=18$
known so far is
$17\,309\,628\,327$ \cite{Aich}.

\bigskip
{\bf \noindent
Acknowledgements:}
The main theorem
was obtained during the \emph{Euroconference on Discrete and Algorithmic
Geometry} held in Anogia, Crete, during August 2000. We thank the organizers of
this meeting, the University of Crete and the Anogia Academic Village for the
invitation and financial support.


\begin{thebibliography}{99}
%
\small
\bibitem{Aich}
\textsc{O.~Aicholzer, H.~Krasser}, 
The point-set order-type database: A collection of
applications and results, in
{\em Proc. $13^{\hbox{th}}$ Canadian Conference on Computational Geometry},
  Waterloo, 2001, pp. 17--20. See also the web page
\textit{Counting triangulations},
{\tt 
\,www.cis.tugraz.at\,/\,igi\,/\,oaich\,/\,triangulations\,/\,counting\,/\,counting.html} 

\bibitem{AHN}
\textsc{O. Aichholzer, F. Hurtado and M. Noy}, On the Number of Triangulations
Every Planar Point Set Must Have, in \textit{Proc.
13th Annual Canadian Conference on Computational Geometry CCCG 2001},
 Waterloo, Canada, 2001, pp. 13--16.

\bibitem{ACNS}
\textsc{M.~Ajtai, V.~Chv\'atal, M.~Newborn and E.~Szemer\'edi},
Crossing-free Subgraphs, \textit{Annals of Discrete Math.},
\textbf{12} (1982), 9--12.
%
\bibitem{BFS}
\textsc{L.~J.~Billera, P.~Filliman and B.~Sturmfels}, 
Constructions and complexity of secondary polytopes, 
\textit{Adv.\ Math.}, \textbf{83} (1990), 155--179.
%
\bibitem{DS}
\textsc{M.~O.~Denny and C.~A.~Sohler},
Encoding a triangulation as a permutation of its point set,
\textit{Proc. 9th Canadian Conf. on Computational Geometry}, (1997), 39--43.
%
\bibitem{GNT}
\textsc{A.~Garc\'{\i}a, M.~Noy and J.~Tejel}, 
Lower bounds on the number of crossing-free subgraphs of $K_N$,
\textit{Comput. Geom.} \textbf{16} (2000), 211--221.
%
%
\bibitem{HN97}
\textsc{F.~Hurtado and M.~Noy}, 
Counting triangulations of almost-convex polygons,
\textit{Ars Combinatoria} \textbf{45} (1997), 169--179.
%
\bibitem{Re}{\sc V. Reiner},
{The generalized Baues problem},
in: {\em New Perspectives in Algebraic Combinatorics}
(L. J. Billera, A. Bj\"orner, C. Greene, R. E. Simion and
R. P. Stanley, eds),
MSRI publications {\bf 38} (1999), Cambridge U. P.,
pp. 293--336.
%
\bibitem{Se}
\textsc{R.~Seidel},
On the number of triangulations of planar point sets,
\textit{Combinatorica}, \textbf{18} (1998), 297--299.
%
%
\bibitem{Sm}
\textsc{W.~D.~Smith}, Studies in Computational Geometry motivated by
Mesh Generation, Ph.~D. Thesis, Princeton Univ. (1989).

\end{thebibliography}
\end{document}